# REJOINDER TO "ANALYSIS OF VARIANCE—WHY IT IS MORE IMPORTANT THAN EVER" BY A. GELMAN

By Andrew Gelman

*Columbia University*

ANOVA is more important than ever because we are fitting models with many parameters, and these parameters can often usefully be structured into batches. The essence of "ANOVA" (as we see it) is to compare the importance of the batches and to provide a framework for efficient estimation of the individual parameters and related summaries such as comparisons and contrasts.

Classical ANOVA is associated with many things, including linear models, F-tests of nested and nonnested interactions, decompositions of sums of squares and hypothesis tests. Our paper focuses on generalizing the assessment, with uncertainty bounds, of the importance of each row of the "ANOVA table." This falls in the general category of data reduction or model summary, and presupposes an existing model (most simply, a linear regression) and an existing batching of coefficients (or more generally "effects," as noted by McCullagh) into named batches.

We thank the discussants for pointing out that more work needs to be done to generalize these ideas beyond classical regression settings with exchangeable batches of parameters. In this rejoinder, we review the essentials of our approach and then address some specific issues raised in the discussions.

**1. General comments.** McCullagh states that we regard "analysis of variance as an algorithm or procedure with a well-defined sequence of computational steps to be performed in fixed sequence." We appreciate this comment, especially in light of Tjur's complaint that it is not clear what our statistical model is. We would like to split the difference and say they are both right: our procedure is indeed performed in a fixed sequence, and the first step is to take a statistical model that must be specified from the outside.

---

A statistical model is usually taken to be summarized by a likelihood, or a likelihood and a prior distribution, but we go an extra step by noting that the parameters of a model are typically batched, and we take this batching as an essential part of the model. If a model is already set up in a fully Bayesian form, our ANOVA step is merely to summarize each batch's standard deviation (whether superpopulation or finite-population; this depends on the substantive context, as discussed by Zaslavsky and in our Section 3.5). If only a likelihood is specified, along with a batching of parameters, we recommend fitting a multilevel model with a variance parameter for each batch, to be estimated from data. Yes, this is an automatic step, and yes, this can be inappropriate in particular cases, but we think it is a big step forward from the current situation in which the analyst must supply redundant information to avoid making inappropriate variance comparisons. Tjur also recognizes our goal of making split-plot and other analyses more understandable to students. More generally, we want to set up a framework where nonstudents can get the correct (classical) answer too (and avoid difficulties such as illustrated in Figure 1)!

Our procedure gives an appropriate answer in a wide range of classical problems, and we find the summary in terms of within-batch standard deviations to be more relevant than the usual ANOVA table of sums of squares, mean squares and F-tests. None of the discussants disputes either of these points, but they all would like to go beyond classical linear models with balanced designs. We provide a more general example in Section 7.2 of our paper (an unbalanced logistic regression) and discuss other generalizations in Section 8.3, but we accept the point that choices remain when implementing ANOVA ideas in nonexchangeable models.

**2. The model comes first.** The discussants raised several important points that we agree with and regret not emphasizing enough in the paper. First, all the discussants, but especially Tjur and McCullagh, emphasize that the model comes first, and the model should ideally be motivated by substantive concerns, not by mathematical convenience and not by the structure of the data or the design of data collection. As noted above, our conception of ANOVA is a way of structuring inferences given that a model has already been fit and that its parameters are already structured into batches. As McCullagh points out, such batches should not necessarily be modeled exchangeably; we defend our paper's focus on exchangeable batches as they are an extremely important special case and starting point (we assume that the coauthor of an influential book on generalized linear models will appreciate the importance of deep understanding of a limited class of models), but note in Section 8.3 that more can be done.

**3. ANOVA is not just for linear models.** Our paper emphasized simple models in order to respond to the unfortunate attitude among many statisticians and econometricians that ANOVA is just a special case of linear regression. Sections 2 and 3 of our paper demonstrate that ANOVA can only be thought of this way if "linear regression" is interpreted to include multilevel models. But ANOVA applies in much more general settings.

The vote-preference example of Section 7.2 is closer to our usual practice, which is to use ANOVA ideas to structure and summarize hierarchical models that have dozens of parameters. In this example, we did not fit a multilevel model because of any philosophical predilections or because we had any particular interest in finite populations, superpopulations or variance parameters. Rather, we sought to capture many different patterns in the data (in particular, state effects to allow separate state estimates, and demographic effects to allow poststratification adjustment for survey nonresponse). The multilevel model allows more accurate inferences—the usual partial pooling or Bayesian rationale [see Park, Gelman and Bafumi (2004)]—and ANOVA is a convenient conceptual framework for understanding and comparing the multiplicity of inferences that result. Compare Figures 6 and 7 to the usual tables of regression coefficients (in this case, with over 50 or 100 parameters) to see the practical advantages of our approach.

Various complications arose naturally in the model fitting stage. For example, state effects started out as exchangeable and then we put in region indicators as state-level predictors. We are currently working on extending these models to time series of opinion polls as classified by states and demographics.

So, even in the example of our paper, the modeling is not as automatic as our paper unfortunately made it to appear. What was automatic was the decision to estimate variance parameters for all batches of parameters and to summarize using the estimated standard deviations. A small contribution, but one that moves us from a tangle of over seventy logistic regression parameters (with potential identifiability problems if parameters are estimated using maximum likelihood or least squares) to a compact and informative display that is a starting point to more focused inferential questions and model improvements.

As the discussants emphasize, in a variety of important application areas we can and should go beyond linear models or even generalized linear models, to include nonlinear, nonadditive and nonexchangeable structures. We have found the method of structuring parameters into batches to be useful in many different sorts of models, including nonlinear differential equations in toxicology, where population variability can be expressed in terms of a distribution of person-level parameters [e.g., Gelman, Bois and Jiang (1996)] and Boolean latent-data models in psychometrics, which have batches of parameters indexed by individuals, situations and psychiatric symptoms [e.g.,

Meulders et al. (2001)]. We cite our own work here to emphasize that we certainly were not trying to suggest that the analysis of variance be restricted to linear models. With modern Bayesian computation, a great deal more is possible, as Hox and Hoijtink point out (and as they have demonstrated in their own applied work). We recommend that practitioners consider ANOVA ideas in summarizing their inferences in these multilevel settings.

**4. ANOVA as a supplement to inferences about quantities of interest.** In a discussion of the example of our Section 2.2.2 (to which we shall return below), Hox and Hoijtink point out that in any specific application an applied researcher will typically be interested in particular treatment effects, or comparisons of treatment effects, rather than in variance components. We agree and thank these discussants for emphasizing this point. As with classical ANOVA, our goal in summarizing variance components is to understand the model as a whole—how important is each source of variation in explaining the data?—as a prelude or accompaniment to more focused inferences. ANOVA may be "more important than ever" but it is intended to add perspective to, not to take the place of, inference for quantities of substantive interest.

To put it another way: if you are already fitting a statistical model, its parameters can probably be grouped into batches, and it is probably interesting to compare the magnitude of the variation of the parameters in each batch. Recent statistical research has revealed many sorts of useful densely parameterized models, including hierarchical regressions, splines, wavelets, mixture models, image models, and so on. However, it can be tricky to understand such models or compare them when they are fit to different datasets. A long list of parameter estimates and standard errors will not necessarily be helpful, partly for simple reasons of graphical display, and partly because an ensemble of point estimates will not capture the variance of an ensemble of parameters [Louis (1984)]. The two examples provided by Zaslavsky illustrate ways in which inferences for variance components can be relevant in applied settings.

Tjur asks about partial confounding and other unbalanced designs. We would simply handle these using Bayesian inference. For example, Section 7.2 gives an example of an unbalanced design. Our paper discussed classical estimates for balanced designs, to connect to classical ANOVA and provide fast calculations for problems like the Internet example, but more generally one can always use full Bayesian computations, as pointed out by Hox and Hoijtink.

**5. Estimation and hypothesis testing.** As Zaslavsky notes, our treatment of ANOVA focuses on estimation of variance components (and, implicitly, of individual coefficients and contrasts via shrinkage estimation), rather than on hypothesis testing. In the application areas in which we have worked, interest has lain in questions of the form, "How important are the effects of factor X?," rather than "Does factor X have an effect?"; see Figure A for an example. (We acknowledge McCullagh's point that our focus is the product of our experiences in environmental, behavioral and social sciences; in other fields, such as genetics, hypotheses of zero effects are arguably more relevant research questions.)

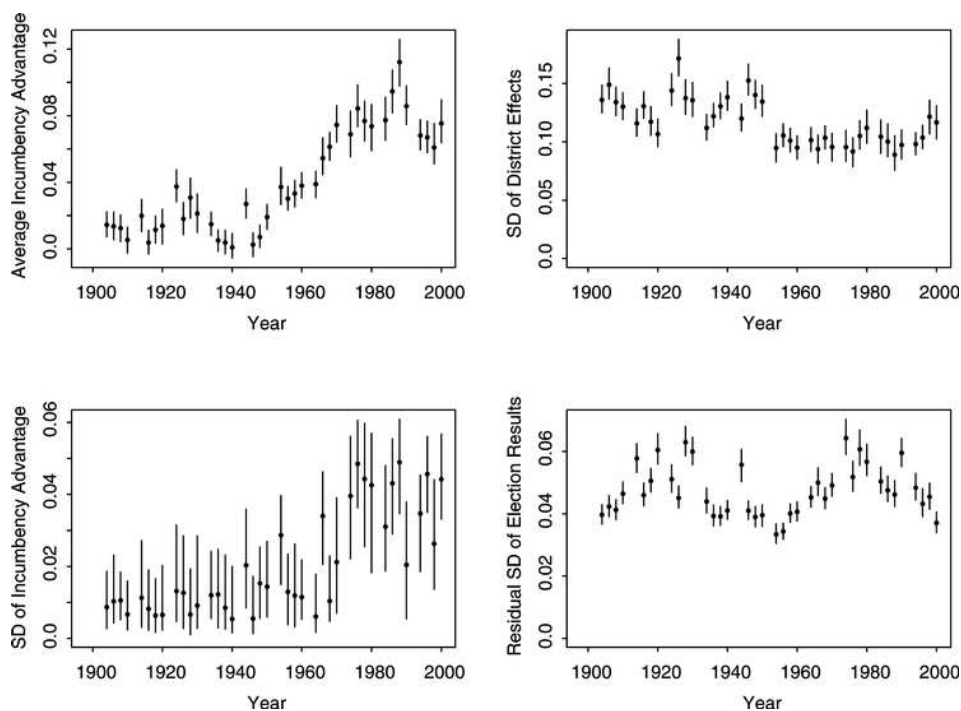

FIG. A. *Estimates and* 95% *intervals for an average treatment effect and three variance parameters for a hierarchical model of elections for the U.S. House of Representatives, fit separately to pairs of successive elections from the past century. The graphs illustrate how we are interested in estimating the magnitude of each source of variation, not simply testing whether effects or variance components equal zero. From Gelman and Huang (2005).*

In settings where hypothesis testing is desired, we agree with Zaslavsky that posterior predictive checking is the best approach, since it allows a hypothesis about any subset of parameters to be tested while accounting for uncertainty in the estimation of the other parameters in the model. Posterior predictive checking can also be applied to the multilevel model as a whole to test assumptions such as additivity, linearity and normality.

**6. Finite-population and superpopulation summaries.** Zaslavsky points out that, in settings where one is interested in generalizing or predicting for new groups, superpopulation summaries are most relevant. We emphasized finite-population summaries in our paper so as to provide more continuity with classical ANOVA. For example, with only five treatment levels, nonzero superpopulation variances are inherently difficult to estimate, a problem that is somewhat ducked by the usual classical analysis which focuses on testing hypotheses of zero variance.

A related issue arises in hierarchical regression models, where the concept of a "contrast" in ANOVA plays the role of a finite-population regression coefficient, while the coefficient in the group-level regression has a superpopulation interpretation. For instance, in the Latin square example shown in Figure 3, suppose we are interested in the linear contrast of treatments A, B, C, D, E. The finite-population contrast is $-2 \cdot \beta_1 + (-1) \cdot \beta_2 + 0 \cdot \beta_3 + 1 \cdot \beta_4 + 2 \cdot \beta_5$, whereas the superpopulation contrast is the appropriately scaled coefficient of $(-2, 1, 0, 1, 2)$ included as a treatment-level predictor in the multilevel model.

A key technical contribution of our paper is to disentangle modeling and inferential summaries. A single multilevel model can yield inference for finite-population and superpopulation inferences. For example, in the example of Section 2.2.2, the structure of the problem implies a model with treatment and machine effects, as noted by Hox and Hoijtink. These authors state that their preferred procedure is "not what [we] had in mind," but they do not fully state what their model is. The key question is: what is the population distribution for the four treatment effect parameters? Our recommendation is to fit a normal distribution with mean and standard deviation as hyperparameters estimated from the data. This is the "superpopulation" standard deviation in the terminology of our Section 3.5; fitting the model would also give inferences for the individual treatment effects and their standard deviation. Hox and Hoijtink question whether the variance of the treatment effects is "an interesting number to estimate"; as we note above in our discussion, we agree with this point but find the general comparison of all the variance parameters to be a useful overall summary (as illustrated by the ANOVA graphs in our paper) without being a replacement for the estimation of individual treatment effects.

To continue with Hox and Hoijtink's discussion of our example: we are not sure what analysis they are suggesting in place of our recommended hierarchical model. One possibility is least-squares estimation for the treatment effect parameters, which would correspond to our hierarchical model with a variance parameter preset to infinity. This seems to us to be inferior to the more general Bayesian approach of treating this variance as a hyperparameter and estimating it from data, and it would also seem to contradict Hox and Hoijtink's opposition to noninformative prior distributions later in their discussion. Another possibility would be a full Bayesian approach with a more informative hyperprior distribution than the uniform distribution that we use. We agree that in the context of any specific problem, a better prior distribution (or, for that matter, a better likelihood) should be available, but we find the normal model useful as a default or starting point.

**7. Fixed and random effects.**  We suspect that statisticians are generally unaware of the many conflicting definitions of the terms "fixed" and "random"; in fact, a reviewer of an earlier version of this paper criticized the multiple definitions in Section 6 as "straw men," which is why we went to the trouble of getting references for each. We are glad that Zaslavsky liked our discussion of fixed and random effects and that McCullagh recognized the "linguistic quagmire."

Hox and Hoijtink would like to define a fixed effect as "a varying effect with components that do not come from the same distribution." This distinction may be important, but we are not hopeful that they will be successful in establishing a new meaning to an already overloaded expression that has at least five other existing interpretations in the statistical literature! Is the phrase "fixed effect" so important that it is worth fighting over this patch of linguistic ground? We use the terms "constant" and "varying" effects because they are unambiguous statements about parameters in a model, and we have the need to communicate with researchers in a wide range of substantive fields. If Hox and Hoijtink find it useful to label sets of effects that are batched but do not come from a common distribution, we recommend they use an unambiguous phrase such as "differently distributed effects" that communicates the concept directly.

Tjur states that our "basic idea seems to be to let all effects enter formally as random effects." We are disappointed to see that he seems to have skipped over Section 6 of our paper! The term "random effect" has no clear (let alone "formal") definition, so we certainly do not consider it to be any part of our basic idea! On the contrary, our basic idea is to recognize that the parameters in a model are not simply a long undifferentiated vector but can be usefully grouped into batches, which in fact are already specified in the classical ANOVA table.

**8. Summary: why is ANOVA important now?**  First, as noted above, if you are already fitting a complicated model, your inferences can be better understood using the structure of that model. We have presented a method for doing so in the context of batches of exchangeable parameters, and we anticipate future developments in other classes of models such as discussed by McCullagh.

Second, if you have a complicated data structure and are trying to set up a model, it can help to use multilevel modeling—not just a simple units-within-groups structure but a more general approach with crossed factors where appropriate. This is the way that researchers in psychology use ANOVA, but they are often ill-served by the classical framework of mean squares and F-tests. We hope that our estimation-oriented approach will allow the powerful tools of Bayesian modeling to be used for the applied goals of inference about large numbers of structured parameters.

## ADDITIONAL REFERENCES

DEPARTMENT OF STATISTICS
COLUMBIA UNIVERSITY
NEW YORK, NEW YORK 10027
USA
E-MAIL: gelman@stat.columbia.edu